\documentclass[dvips,preprint]{imsart}

\RequirePackage[OT1]{fontenc}
\RequirePackage{amsthm,amsmath}
\RequirePackage[colorlinks,citecolor=blue,urlcolor=blue]{hyperref}
\usepackage{amsfonts, amsmath, bbm, amssymb, latexsym, graphics, graphicx}

\startlocaldefs
\numberwithin{equation}{section}
\theoremstyle{plain}
\newtheorem{theorem}{Theorem}[section]
\newtheorem{lemma}[theorem]{Lemma}
\newtheorem{proposition}[theorem]{Proposition}
\newtheorem{corollary}[theorem]{Corollary}
\newtheorem{definition}[theorem]{Definition}

\endlocaldefs

\DeclareMathOperator{\C}{C}
\DeclareMathOperator{\CV}{CV}

\DeclareMathOperator{\vol}{Vol}
\DeclareMathOperator{\tube}{Tube}

\DeclareMathOperator{\sign}{sign}

\def\N{\mathbb{N}}

\def\R{\mathbb{R}}
\def\Z{\mathbb{Z}}

\def\N{\mathcal{N}}

\def\e{\varepsilon}

\def\l{\lambda}

\def\w{\omega}

\def\vp{\varphi}

\def\ldc{\lfloor d\chi \rfloor}
\def\udc{\lceil d\chi \rceil}

\newcommand{\E}{\mathbb{E}} 

\newcommand{\mean}[1] {\E\left\{{#1}\right\}}

\newcommand{\prob}[1]{\mathbb{P}\left(#1\right)}
\newcommand{\ind}{\boldsymbol{\mathbbm{1}}} 

\newcommand{\LL}{{\mathcal{L}}}
\newcommand{\MM}{{\mathcal{M}}}

\newcommand{\set}[1]{\left\{#1\right\}}
\newcommand{\sbrk}[1]{\left[#1\right]}

\newcommand{\norm}[1]{\left\|#1\right\|}
\newcommand{\param}[1]{\left(#1\right)}
\newcommand{\abs}[1] {\left| {#1}\right|}
\newcommand{\floor}[1] {\left\lfloor{#1}\right\rfloor}
\newcommand{\ceil}[1] {\left\lceil{#1}\right\rceil}
\newcommand{\iprod}[1] {\left< {#1}\right>}

\newcommand{\EC}{{Euler characteristic}}
\newcommand{\EI}{{Euler integral}}

\newcommand{\comb}[2]{{{#1}\atopwithdelims[]{#2}}}
\usepackage{xcolor}

\begin{document}

\begin{frontmatter}
\title{Euler Integration of Gaussian Random Fields and Persistent Homology}
\runtitle{Euler Integration of Gaussian Random Fields}

\begin{aug}
\author{\fnms{Omer} \snm{Bobrowski}\thanksref{t1,t2,t3}\ead[label=e1]{bober@tx.technion.ac.il}},
\author{\fnms{Matthew Strom} \snm{Borman}\thanksref{t1}\ead[label=e2]{borman@math.uchicago.edu}}

\thankstext{t1}{Research supported in part by National Sceince Foundation, DMS-0852227.}
\thankstext{t2}{Research supported in part by US-Israel Binational Science Foundation, 2008262.}
\thankstext{t3}{Supported in part by the Adams Fellowship Program of the Israel Academy of Sciences and Humanities.}
\runauthor{Bobrowski and Borman}

\affiliation{Technion\thanksmark{m1} and University of Chicago\thanksmark{m2}}

\address{Omer Bobrowski\\
Department of Electrical Engineering\\
Technion, Israel Institute of Technology \\
Haifa, 32000\\
Israel\\
\printead{e1}\\
\phantom{E-mail:\ }}

\address{Matthew Strom Borman\\
Department of Mathematics\\
University of Chicago\\
Chicago, Illinois 60637\\
USA\\
\printead{e2}\\
\phantom{E-mail:\ }}
\end{aug}

\begin{abstract}
In this paper we extend the notion of the Euler characteristic to persistent homology and give the relationship between the Euler integral of a function and the Euler characteristic of the function's persistent homology.  We then proceed to compute the expected Euler integral of a Gaussian random field using the Gaussian kinematic formula and obtain a simple closed form expression.  This results in  the first explicitly computable mean of a quantitative descriptor for the persistent homology of a Gaussian random field.
\end{abstract}

\begin{keyword}[class=AMS]
\kwd{60G15, 55N35}
\end{keyword}

\begin{keyword}
\kwd{Persistent homology, barcodes, Betti numbers, Euler characteristic, random fields, Gaussian processes, Gaussian kinematic formula.}
\end{keyword}

\end{frontmatter}

\section{Introduction}

The Euler characteristic, denoted by $\chi(X)$, is an invariant of nice topological spaces $X$, which is additive in the sense that there is inclusion-exclusion: \[\chi(A \cup B) = \chi(A) + \chi(B) - \chi(A \cap B).\]
Since this is the heart of a measure, one can make sense of integration with respect to the
Euler characteristic.  This new integration theory has many of the properties of Lebesgue integration,
such as a Fubini theorem and integral transforms (e.g.\! Fourier and Radon transforms \cite{viro1988some}).
Aside from the theoretical aspects, the Euler integral has recently started to give rise to applications. For example, \cite{baryshnikov2009target} presents a target enumeration scheme which employs the {\EI} to perform data aggregation in sensor networks.

One major drawback of the original Euler integral is that it was only defined for integer valued functions.
In \cite{baryshnikov2009euler}, two extensions were given, called the \textit{upper} and \textit{lower} Euler integrals, which are defined for a large class of real valued functions. In this paper we will use the upper Euler integral, but everything we do has a lower integral analogue.
This integral operator, while having some drawbacks (particularly a lack of linearity), possesses a compelling Morse theoretic interpretation (see \cite{baryshnikov2009euler}).
In this paper we extend this Morse theoretic interpretation to general tame functions.
Additional information on the stratified Morse theory used in this paper can be found in \cite{Goresky:Macpherson}.  For algebraic topology preliminaries on the Euler characteristic and homology see \cite{hatcher2002algebraic}.

We will use recent developments in the theory of Gaussian random fields to compute the expected value of the upper Euler integral of Gaussian and Gaussian related random fields.
Random fields are stochastic processes defined over a topological space $X$.
 Gaussian fields are characterized by having finite dimensional distributions which are multivariate normal.  A field generated by applying a function to a vector valued Gaussian field is called a \textit{Gaussian related field}.
Let $f$ be a $k$-dimensional, Gaussian field on a nice space $M$ of dimension $d$, $G:\R^k\to\R$ a piecewise $C^2$ function, $g=G\circ f$. Assume that $f$ has zero mean and constant unit variance. Then, under mild regularity conditions, our most general result states that
\begin{equation}\label{eq:general_intro}
\mean{\int_M{g \udc}} = \chi(M)\mean{g} - \sum_{j=1}^d{(2\pi)^{-j/2}\LL_j(M)\int_{\R}{\MM_j(D_u) du}},
\end{equation}
where $D_u = G^{-1}(-\infty,u]$, $\mean{g}:= \mean{g(t)}$ (for any $t\in M$), and the $\LL_j,\MM_j$ are geometric characteristics known as Lipschitz-Killing curvatures and Gaussian Minkowski functionals, respectively.

There is a strong link between the Euler integral and persistent homology. Briefly, the persistent homology of a real valued function $f$ tracks changes in the homology of sublevel sets $f^{-1}(-\infty,u]$. As the level $u$ changes, new homology elements (i.e.\! ``holes") are born and others die. Persistent homology keeps a record of this birth-death process.
The theory of persistent homology was created as a way to describe the topology of data, but to date most results about persistent homology have been algebraic in nature as opposed to probabilistic or statistical.  For the theory of persistent homology to develop into a powerful applied tool there will need to be theorems about persistent homology in the random setting.  A first step in this direction was \cite{bubenik2007statistical}, where the expected barcode of the persistent homology for points sampled from a circle was computed.  We will relate persistent homology to the Euler integral via a parameter we call \textit{the Euler Characteristic of the Persistent Homology}. Thus, we obtain the first known general tool that allows one to make probabilistic statements about the persistent homology generated by sublevel sets.  Furthermore, our computation of the expected Euler integral for Gaussian random fields provides the first quantitative measure of the persistent homology of random fields in the multi-dimensional setting.

Another significant and surprising corollary of the formula \eqref{eq:general_intro} is related to the signed sum of critical values of Gaussian random fields. If $f:M \to \R$ is a Gaussian random field and $M$ is a closed manifold
we prove that
\[
\mean{\sum_{p \in C(f)} (-1)^{\mu(p)} f(p)} = -\frac{\LL_1(M)}{\sqrt{2 \pi}},
\]
where $\C(f)$ is the set of critical points of the field $f$, and $\mu(p)$ is the Morse index of $f$ at $p$.
In other words, the expected signed sum of critical values of a Gaussian field does not scale according to the volume of the space as one might expect, but rather according to a one dimensional measure of the space.

Finally, using \eqref{eq:general_intro} we can extend the application discussed in \cite{baryshnikov2009target} to the situation where the sensor field is contaminated by additive Gaussian noise. In this particular case the integral remains additive, which enables us to suggest a noise reduction scheme.

\section{The Euler Integral}

The {\EC} is an additive operator on compact sets. Therefore, it is tempting to consider $\chi$ as a measure and integrate with respect to it. The main problem in doing so is that $\chi$ is only finitely additive.
At first (see \cite{viro1988some}), integration with respect to the {\EC} was defined for a small set of functions called \textit{constructible functions} defined by
\[
C\!F(X) = \set{\left.f(x) = \sum_{k=1}^n a_k\, \ind_{A_k}(x)\,\right|\, a_k\in \Z,\, A_k\, \textrm{disjoint tame subsets of }X}
\]
where `tame' means having a finite \EC.
For this set of functions we can define the {\EI} analogously to the Lebesgue integral. Let $f(x) = \sum_{k=1}^n a_k \ind_{A_k}(x)$, and define
\[
\int_X f d\chi = \sum_{k=1}^n a_k \chi(A_k).
\]
This integral has many nice properties, similarly to those of the Lebesgue integral, such as linearity and a version of the Fubini theorem (see \cite{viro1988some}).
However, as mentioned above, the {\EC} is not countably additive, and therefore we cannot continue from here by approximating other functions using functions in $C\!F(X)$.

In \cite{baryshnikov2009euler} two extensions were suggested for the {\EI} of real valued functions using the notion of a definable function over an $\mathcal{O}$-minimal structure (see \cite{baryshnikov2009euler} and references therein for more background). Let $\floor{x}$ and $\ceil{x}$ be the floor and ceiling values of $x$ respectively.  In the $\mathcal{O}$-minimal language, if $f: X \to \R$ is a definable function on a definable space $X$, then an important property is that both $\floor{f}$ and
$\ceil{f}$ are constructible functions and hence have well defined Euler integrals.  This leads to the following Riemann-sum like definition:
\begin{definition}[\cite{baryshnikov2009euler}]\label{def:BG def}
The lower {\EI} of a definable function $f: X \to \R$ on a definable space $X$ is defined by
\begin{equation} \label{eq:bg_lei}
\int_X f \ldc = \lim_{n\to\infty} n^{-1} \int_X \floor{n f} d\chi,
\end{equation}
and the upper {\EI} is defined by
\begin{equation} \label{eq:bg_uei}
\int_X f \udc = \lim_{n\to\infty} n^{-1} \int_X \ceil{n f} d\chi.
\end{equation}
\end{definition}
These two extensions coincide with the original {\EI} for constructible functions.  For other functions they might be completely different.

Unfortunately, it is not clear if Gaussian random fields $f: X \to \R$ can be made to fit inside an
$\mathcal{O}$-minimal setting.  Therefore in this work we will use the following simplified definition of a \textit{tame} function.
\begin{definition}\label{def:BB tame}
	A continuous function $f: X \to \R$ on a compact topological space $X$
	with a finite {\EC} is \emph{tame} if the homotopy types of
	$f^{-1}(-\infty, u]$ and $f^{-1}[u, \infty)$ change only finitely many times as $u$ varies over $\R$
	and the Euler characteristic of each set is always finite.
\end{definition}
We can now take our definition of the lower and upper Euler integrals to be as follows.
\begin{definition}\label{prop:ei_computation}
	If $f: X \to \R$ is a tame function, then the lower and upper Euler integrals are defined by
    \begin{align}\label{eq:lei}
		\int_X f \ldc &= \int_{u=0}^\infty \sbrk{\chi(f \geq u) - \chi(f < -u)}\, du \\
        \label{eq:uei}
        \int_X f \udc &= \int_{u=0}^\infty\sbrk{ \chi(f > u) - \chi(f \le -u)}\, du,
    \end{align}
where $\chi(f \ge u) := \chi\param{f^{-1}[u,\infty)},\ \chi(f < u) := \chi(X)- \chi(f \ge u)$, etc.
\end{definition}
These equations appear as Proposition 2 in \cite{baryshnikov2009euler}, where they were proven using Definition \ref{def:BG def}.  Thus, Definition~\ref{def:BG def} and Definition~\ref{prop:ei_computation}
agree when applied to a function that is both definable in the $\mathcal{O}$-minimal sense and tame in the sense of Definition~\ref{def:BB tame}.  Rather than viewing \eqref{eq:lei} and \eqref{eq:uei} as consequences of the $\mathcal{O}$-minimal definitions \eqref{eq:bg_lei} and \eqref{eq:bg_uei},
we will take \eqref{eq:lei} and \eqref{eq:uei} to be the definition of the Euler integrals for functions that are tame in the sense of Definition~\ref{def:BB tame}.

\subsection{The Euler Integral and Morse Theory}
In \cite{baryshnikov2009euler} the Euler integral was given a stratified Morse theory interpretation.  A corollary of this approach was that if $f: M \to \R$ is a Morse function, and $M$ is a closed manifold, then
\begin{equation}\label{eq:ei_critical_points}
	\int_M f \udc = \sum_{p \in \C(f)} (-1)^{\mu(p)} f(p),
\end{equation}
where $\C(f)$ is the set of critical points of $f$ and $\mu(p)$ is the index of $p$ as a critical point.
In our language of tame functions, we have the following proposition.
\begin{proposition}\label{prop:crit_val}
	Let $f: X \to \R$ be a tame function and let $\CV(f)$
	be the set of values where the homotopy type of $f^{-1}(-\infty, u]$
	changes (the critical values of $f$).
	Then
	\[
		\int_X f \udc = \sum_{v\in \CV(f)} \Delta_\chi(f,v)\, v,
	\]
	where $\Delta_\chi(f,v) = \chi(f \leq v+\e) - \chi(f \leq v-\e)$, for sufficiently small $\e$,
	is the change in the Euler characteristic of $f^{-1}(-\infty, u]$ as $u$ passes through the critical
	value $v$.
\end{proposition}
\begin{proof}
    Label the critical values $\CV(f) = \set{v_1,\dots,v_n}$ in increasing order such that
    $v_1 < \cdots < v_i < 0 \le v_{i+1} < \cdots < v_n$. If $v_k < u < v_{k+1}$, then via a telescoping sum
\begin{align*}
		\chi(f \le u) &= \Delta_\chi(f,v_1) + \cdots + \Delta_\chi(f,v_k),\\
		\chi(f > u) &= \chi(M) - \chi(f \le u) = \Delta_\chi(f,v_{k+1}) + \cdots + \Delta_\chi(f,v_n).
\end{align*}
	Therefore for $u \in [0, \infty)$ and $u \not= \pm v_j$,
\begin{align*}
		\chi(f > u) &= \sum_{j=i+1}^n \Delta_\chi(f,v_j)\, \ind_{[0,v_j]}(u) \quad \mbox{and}\\
		\chi(f \le -u) &= \sum_{j=1}^i \Delta_\chi(f,v_j)\, \ind_{[0, -v_j]}(u).
\end{align*}
	Thus
	\[
		\int_M f \udc = \int_{u=0}^\infty \param{\chi(f > u) - \chi(f \le -u)}du
			= \sum_{j=1}^n v_j\, \Delta_\chi(f,v_j),
	\]
	as desired.
\end{proof}

This recovers the Morse theoretic viewpoint, since if $f: M \to \R$ is a Morse function then
Morse theory says that the {\EC} changes by the addition of $(-1)^k$ as $f(-\infty, u]$
passes through a critical point of index $k$.

\begin{corollary}\label{cor:alter_sum_Ga}
Let $f: X \to \R$ be tame, satisfying the conditions of Proposition~\ref{prop:crit_val}.
If $a$ is not a critical value of $f$, then
\[
	\sum_{\stackrel{v\in \CV(f)}{v< a}} \Delta_\chi(f, v) v
		= \int_X G_a(f) \udc + a\,\chi(f \leq a) - a\,\chi(X)
\]
where $G_a(x) = \min(x,a)$.

\end{corollary}
\begin{proof}
	If $v_1 < \cdots < v_n$ are the critical values of $f$ and
	$a$ is such that $v_k < a < v_{k+1}$, then $G_a(f)$ has critical values at
	$ v_1 < \cdots < v_k < a$.
	By Proposition~\ref{prop:crit_val},
	\begin{align*}
		\int_X G_a(f) \udc &= \sum_{j=1}^k v_j\,\Delta_\chi(G_a(f), v_j) + a\,\Delta_\chi(G_a(f), a)\\
		&= \sum_{j=1}^k v_j\,\Delta_\chi(f, v_j) + a(\chi(M)-\chi(f \leq a)).
	\end{align*}
	This gives us the desired result.
\end{proof}

\section{Gaussian Random Fields and the Gaussian Kinematic Formula}
%

There has been extensive effort over the past few years to study the sample paths
of smooth random fields from a general Riemannian manifold $M$ to $\R^k$.
Specific examples in which this approach has had practical importance occur when $M$ is a 3-dimensional brain or a 2-dimensional cortical surface. The basic (random) geometrical objects studied were the excursion sets of the random fields, namely $f^{-1}(D)$, for nice subsets $D$ of $\R^k$, and the tools for quantifying these sets were those of differential geometry.
The theory of this subject has developed rapidly over the past few years (see \cite{adler2007random,taylor2001thesis,taylor2006gaussian}).
One of its most powerful results is an explicit expression for the mean value of all Lipschitz-Killing curvatures (among them is the \EC) of excursion sets for centered (i.e.\! $\mean{f(t)}=0$), constant variance, $C^2$, Gaussian random fields. The result presented in \cite{adler2007random} links random field theory with integral and differential geometry, and leads to approximations of important objects in
probability and statistics, such as the exceedence probabilities $\prob{\sup_M f(t) > u}$  \cite{adler2000excursion, taylor2005validity}. The main theorem in \cite{adler2007random} is called the \textit{Gaussian kinematic formula} (GKF). We state it here without getting fully involved with the details.

\begin{theorem}[The GKF, {\cite[Theorem 15.9.4]{adler2007random}}]\label{thm:gkf}
Let $M$ and $D\subset \R^k$ be regular stratified spaces with $M$ compact and $D$ closed.
Let $f=(f_1,\ldots,f_k):M\to\R^k$ be a $k$-dimensional Gaussian field, with $iid$ components all having zero mean, unit variance and such that with probability one $f_j$ is a stratified Morse function. Then
\[
\mean{\LL_i(f^{-1}(D))} = \sum_{j=0}^{\dim M-i} \comb{i+j}{j}(2\pi)^{-j/2}\LL_{i+j}(M)\MM_j(D),
\]
where:
\begin{itemize}
\item $\LL_i(\cdot)$ is the $i$-th Lipschitz-Killing curvature of $X$ (with respect to the metric defined in \eqref{eq:metric} below),
\item $\MM_i(\cdot)$ is the $i$-th Gaussian Minkowski functional,
\item $\comb{n}{k} = \binom{n}{k} \frac{\w_n}{\w_k \w_{n-k}}$, and  $\w_n$ is the volume of the $n$-dimensional unit ball.
\end{itemize}
\end{theorem}

The definitions of regular stratified spaces, the Lipschitz-Killing curvatures, and the Gaussian Minkowski functionals are given in \cite[Definition 9.2.2, Definition 10.7.2, Corollary 10.9.6]{adler2007random}.
While the precise definition is involved, examples of regular stratified spaces are
closed manifolds, compact manifolds with boundary, and products of regular stratified spaces.
Conditions are given in \cite[Corollary 11.3.5]{adler2007random} on the covariance function of $f$, such that sample paths are stratified Morse functions with probability one (and so in particular are tame functions).
See \cite{adler2007random} as well for a reference for most of the following facts.

Lipschitz-Killing curvatures are
geometric objects that depend on the choice of a Riemannian metric on $M$, such that $\LL_k(M)$ is a measure of the $k$-dimensional `size' of $M$.   This means that if we scale the metric by $\l$, then $\LL_k(M)$ scales by $\l^k$. For a large class of spaces, which include smooth manifolds and convex compact regions, if $X \subset \R^n$ is given the Euclidean metric, then the following tube formula holds for sufficiently small $r$
\begin{equation}\label{eq:tube_formula}
	\mu(B(X, r)) = \sum_{j=0}^n \w_{j} \LL_{n-j}(X) r^{j},
\end{equation}
where $\mu(B(X,r))$ is the Lebesgue measure of the tube of radius $r$ about $X$.  It turns out that the $\LL_j(X)$ are independent of how $X$ is isometrically embedded into $\R^n$ (and of $n$).
For a $d$-dimensional space, $\LL_d(X)$ is its Riemannian volume and $\LL_0(X)$ is always its Euler characteristic.
Here are few examples, if $X \subset \R^2$ is convex and compact then
\[
	\LL_0(X) = 1, \quad \LL_1(X) = \mbox{(perimeter of $X$)}/2,\quad \LL_2(X) = \textrm{area}(X).
\]
Now suppose that $M$ is a closed $d$-dimensional manifold.
If $M$ is odd dimensional then the even $\LL_k$ vanish, while if $M$ is even dimensional then the odd $\LL_k$ vanish.  When the parities match, $\LL_{d-j}(M)$ for even $j$ is given by an integral of an expression involving the curvature tensor.  In the case that $M$ has constant sectional curvature $\kappa$, then
\[
	\LL_{d-j}(M) = (4\pi)^{-j/2} \frac{d!}{(d-j)!(j/2)!}\, \kappa^{j/2}\, \vol(M).
\]

The Gaussian Minkowski functionals ($\MM_i(D)$) are also defined using a tube formula similar to the one in \eqref{eq:tube_formula}, but using the Gaussian measure on $\R^k$ instead of Lebesgue measure.

If $f: M \to \R$ is a mean zero Gaussian random field, then it is completely determined by
its covariance function $C: M \times M \to \R$, where $C(s,t) = \mean{f(s)f(t)}$, and we can use this to define a metric $g$ on $M$ by
\begin{equation}\label{eq:metric}
	g_t(X_t, Y_t) = \mean{X_t(f)Y_t(f)} = X_sY_tC(s,t)\mid_{s=t},
\end{equation}
where $X$ and $Y$ are vector fields on $M$ defined near $t \in M$.
It is with respect to this metric that the Lipschitz-Killing curvatures are computed in the GKF,
and also the metric with respect to which we require $M$ to be bounded.  It turns out that any Riemannian metric can be realized as coming from a variance one Gaussian random field.

We will be particularly interested in the special case of the GKF for which $i=0$.  Since $\LL_0$ is just the {\EC} $\chi$, the GKF implies that
\[
\mean{\chi(f^{-1}(D))} = \sum_{j=0}^{\dim M} (2\pi)^{-j/2}\LL_j(M) \MM_j(D).
\]

\section{The Euler Integral of Gaussian Random Field} \label{sec:ei_grf}

Let $M$ be a stratified space and let $g: M \to \R$ be a Gaussian or Gaussian related random field.
We are interested in computing the expected value of the Euler integral of the field $g$ over $M$.
While we focus on the upper Euler integral, everything we do has a lower Euler integral analogue.
The following result is a corollary of the GKF and Proposition~\ref{prop:ei_computation}.
\begin{theorem}\label{thm:mean_ei_general}
Let $M$ be a compact $d$-dimensional stratified space, and let $f:M\to\R^k$ be a $k$-dimensional Gaussian random field satisfying the GKF conditions. For a piecewise $C^2$ function $G:\R^k \to \R$,
let $g=G\circ f$. Setting
$D_u = G^{-1}(-\infty,u]$, we have
\begin{equation}\label{eq:mean_ei}
\mean{\int_M{g\udc}} = \chi(M)\mean{g} - \sum_{j=1}^d{(2\pi)^{-j/2}\LL_j(M)\int_{\R}{\MM_j(D_u) du}}
\end{equation}
where $\mean{g} = \mean{g(t)}$ ($g(t)$ has a constant mean).
\end{theorem}

The difficulty in evaluating the expression above lies in computing the Minkowski functionals
$\MM_j(D_u)$.
In Sections \ref{sec:mean_ei_real} and \ref{sec:mean_ei_vec} we present a few cases where they have been computed, which allows us to simplify \eqref{eq:mean_ei}.

\begin{proof}
By Proposition \ref{prop:ei_computation}
\begin{align*}
\int_M g \udc &= \int_{0}^\infty \param{\chi(g > u) - \chi(g \le -u)}du \\
&= \int_{0}^\infty \param{\chi(M)-\chi(g \le u)}du - \int_{-\infty}^0\chi(g \le u)du .
\end{align*}
Therefore,
\begin{equation}\label{eq:upper_ei_sublevels}
\begin{split}
\mean{\int_M g \udc} &= \int_{0}^\infty \param{\chi(M)-\mean{\chi(g \le u)}}du\\
&\mbox{\quad\quad\quad}- \int_{-\infty}^0\mean{\chi(g \le u)}du.
\end{split}
\end{equation}
Replacing $D$ with $D_u$ in the GKF (Theorem \ref{thm:gkf}) yields
\[
\mean{\chi(g \le u)} = \mean{\chi(f^{-1}(D_u))}
=\sum_{j=0}^{d}(2\pi)^{-j/2}\LL_{j}(M)\MM_j(D_u).
\]
Substituting this formula into \eqref{eq:upper_ei_sublevels} yields,
\[
\mean{\int_M{g \udc}} = \sum_{j=0}^d{(2\pi)^{-j/2}\rho_j \LL_j(M)},
\]
where
\[
\rho_j = \begin{cases}
 -\int_{\R}{\MM_j(D_u) du} & j > 0,\\
 \int_0^\infty{\param{1-\MM_0(D_u)}du}-\int_{-\infty}^0{\param{\MM_0(D_u) du}} & j=0.
\end{cases}
\]
The expression for $\rho_0$ can be further simplified. Let $X$ be a standard normal variable and $Y=G(X)$, then
\[
\MM_0\param{D_u} = \gamma_k(D_u) =  \prob{X \in D_u} = \prob{Y\le u}.
\]
Therefore,
\begin{align*}
\rho_0 &= \int_0^\infty {\param{1-\prob{Y \le u} }du} - \int_{-\infty}^0 {\prob{Y \le u}du}\\
&= \int_0^\infty {\prob{Y > u} du} - \int_{-\infty}^0 {\prob{Y \le u}du}\\
&= \mean{Y}.
\end{align*}
Since for all $t$, $f(t)\sim \N(0,1)$, we can replace $Y$ with $G(f(t)) = g(t)$.
Finally, recalling that $\LL_0 \equiv \chi$ completes the proof.
\end{proof}

\subsection{Real Valued Fields}\label{sec:mean_ei_real}

For real valued fields we can improve Theorem~\ref{thm:mean_ei_general} by
computing the terms $\MM_j(D_u)$ that appear in \eqref{eq:mean_ei}.  First, we need to recall some facts about the family of Hermite polynomials. For $n \geq0$, the $n$-th Hermite polynomial is defined as
\[
	H_n(x) = (-1)^n \varphi(x)^{-1} \frac{d^n}{dx^n} \varphi(x),
\]
where
$\varphi(x) = (2\pi)^{-1/2}e^{-x^2/2}$ is the density of the standard Gaussian distribution. This family of polynomials is orthogonal under the inner product on functions $f,g:\R \to \R$
\[
	\iprod{f,g} = \int_{\R}{f(x)g(x)\vp(x)}dx.
\]
A useful convention is 
\[
	H_{-1}(x) = \varphi(x)^{-1}\int_x^\infty \varphi(u) du.
\]

\begin{theorem}\label{thm:mean_ei_real}
Let $M$ be a compact $d$-dimensional stratified space, and let $f:M\to\R$ be a real valued Gaussian random field satisfying the GKF conditions. Let $G: \R \to \R$ be piecewise $C^2$ and $g=G \circ f$. Then
\[
\mean{\int_M{g \udc}} = \chi(M)\mean{g} + \sum_{j=1}^d{(-1)^j\LL_j(M)} \frac{\iprod{H_{j-1}, (\sign(G'))^j G'}}{(2\pi)^{j/2}}.
\]

\end{theorem}

In the case that the function $G$ is strictly monotone, this can be simplified by
using the fact that $\sign(G')$ is constant and then integrating by parts.
\begin{corollary}\label{cor:mean_ei_real_monotone}
Let $f$ be as in Theorem \ref{thm:mean_ei_real}, and $G$ be a strictly increasing function. Then
\[
\mean{\int_M{g \udc}} = \sum_{j=0}^d {(-1)^j\LL_j(M) \frac{\iprod{H_j,G}}{(2\pi)^{j/2}}}.
\]
If $G$ is strictly decreasing then,
\[
\mean{\int_M{g \udc}} = \sum_{j=0}^d {\LL_j(M) \frac{\iprod{H_j,G}}{(2\pi)^{j/2}}}.
\]
\end{corollary}

To prove Theorem \ref{thm:mean_ei_real} we will need the following calculus lemma, which is a special case of Federer's coarea formula.
\begin{lemma} \label{lem:change_of_variables}
Let $h:\R\to\R$ be an integrable function and let $G:\R\to\R$ be a piecewise
differentiable continuous function that is nondifferentiable on a discrete set.
Then
\[
\int_{\R}{h(x)\abs{G'(x)}dx} = \int_{\R} { \left(\sum_{x \in G^{-1}(t)} {h(x)}\right)dt}.
\]
\end{lemma}

\begin{proof}(Theorem \ref{thm:mean_ei_real})
By Theorem~\ref{thm:mean_ei_general}, it suffices to show that
\begin{equation}\label{eq:proof 8}
\int_{\R}{\MM_j(D_u) du} = (-1)^j\iprod{H_{j-1}, (\sign(G'))^j G'},
\end{equation}
for $j \geq 1$, where $D_u = G^{-1}(-\infty, u]$.
Since $G$ is continuous, we can write the inverse image of $(-\infty,u]$ as a disjoint union of closed intervals
\[
D_u = \bigcup_i [a_i,b_i]
\]
where we allow one $a_i$ to be $-\infty$ and one $b_i$ to be $\infty$. Note that for all the finite values we have $G(a_i) = G(b_i) = u$, $G'(a_i)<0$ and $G'(b_i)>0$.

For small enough $\rho$ we have
\[
\tube \param{D_u,\rho} = \bigcup_i [a_i-\rho,b_i+\rho].
\]
Therefore
\begin{equation}\label{eq:tube}
\gamma_k\param{\tube \param{D_u,\rho}} = \sum_i \param{\Phi(b_i+\rho) - \Phi(a_i - \rho)},
\end{equation}
where $\Phi(x) = \int_{-\infty}^x \varphi(u)du$.
The Taylor expansion of $\Phi(x+\rho)$ in $\rho$ is
\begin{equation}\label{eq:hermite}
\Phi(x+\rho) = \Phi(x) + \sum_{j=1}^\infty {\frac{\rho^j}{j!}(-1)^{j-1}H_{j-1}(x) \varphi(x)},
\end{equation}
so in particular $\MM_j(-\infty, x] = (-1)^{j-1}H_{j-1}(x) \varphi(x)$.
Therefore we conclude that for $j \ge 1$
\begin{equation}\label{eq:minkowski_j}
\MM_j\param{D_u} = \sum_i \param{(-1)^{j-1}H_{j-1}(b_i)\varphi(b_i) + H_{j-1}(a_i)\varphi(a_i)}.
\end{equation}
Note that if $b_i=\infty$ (or $a_i=-\infty$) its contribution to the volume of the tube in \eqref{eq:tube}
is independent of $\rho$ (1 or 0 respectively). Thus, it will affect only $\MM_0$ and we can assume all the $a_i$ and $b_i$ in \eqref{eq:minkowski_j} are finite and hence $\bigcup_i \{a_i, b_i\} = G^{-1}(u)$.

If $j$ is odd, then from \eqref{eq:minkowski_j} we have that
\[
\MM_j\param{D_u} = \sum_i \param{H_{j-1}(b_i)\varphi(b_i) + H_{j-1}(a_i)\varphi(a_i)} = \sum_{x\in G^{-1}(u)}{H_{j-1}(x)\varphi(x)}.
\]
If $j$ is even, then
\begin{align*}
\MM_j\param{D_u} &= \sum_i\param{-H_{j-1}(b_i)\varphi(b_i) + H_{j-1}(a_i)\varphi(a_i)}\\
 &= -\sum_{x\in G^{-1}(u)}{\sign(G'(x))H_{j-1}(x)\varphi(x)}.
\end{align*}

For the case that $j$ is odd, apply Lemma \ref{lem:change_of_variables} to get
\begin{align*}
\int_{\R} {\MM_j\param{D_u}du} &= \int_{\R} {\left(\sum_{x\in G^{-1}(u)}{H_{j-1}(x)\varphi(x)}\right)du}\\
& = \int_{\R}{H_{j-1}(x)\varphi(x)|G'(x)|dx} \\
&= \iprod{H_{j-1},|G'|} \\
&= (-1)^{j-1} \iprod{H_{j-1},(\sign(G'))^j G'}.
\end{align*}
So we have proved \eqref{eq:proof 8}, when $j$ is odd.
If $j$ is even, a similar calculation gives the desired result.
\end{proof}

\subsection{Vector Valued Fields}\label{sec:mean_ei_vec}
When $f$ is a vector valued Gaussian field, it can be difficult to evaluate the Minkowski functionals
$\MM_j$.  In two cases though, it is possible to compute the mean \EI.

\subsubsection{The $\chi^2$ case}
Let $M$ be a compact $d$-dimensional manifold.  A $\chi^2$ field with $k$ degrees of freedom is of the form $g=G\circ f$, where
$f=(f_1,\ldots,f_k): M \to \R^k$ is a Gaussian random field with $iid$ components,  and $G(x_1,\ldots,x_k) = \sum_{i=1}^k{x_i^2}$.
\begin{theorem}\label{thm:mean_ei_chi}
The mean Euler integral for a $\chi^2$ field with $k$ degrees of freedom, with
$k\ge d$, is given by
\[
\mean{\int_M{g \udc}} = k\,\LL_0(M) - \frac{2}{\sqrt{\pi}}\,\frac{\Gamma(\frac{k+1}{2})} {\Gamma(\frac{k}{2})}\LL_1(M) + \frac{1}{\pi} \LL_2(M).
\]
\end{theorem}
\begin{proof}
First note that in this case, $\MM_j\param{D_u} = \MM_j\param{G^{-1}(-\infty, u]} = 0$
when $u < 0$ since $G$ is nonnegative.
In \cite[Section 15.10.2]{adler2007random} it is shown that for $k \ge d$ and $j \geq 1$
\[
\MM_j\param{D_u}= \left.\frac{d^{j-1} p_k(x)}{dx^{j -1}} \right|_{x=\sqrt{u}}
\quad \mbox{where} \quad p_k(x) = \frac{x^{k-1} e^{-x^2/2}}{\Gamma(k/2) 2^{(k-2)/2}}.
\]
Therefore,
\begin{align*}
\int_\R {\MM_j(D_u)du} =  \int_0^\infty{ \left.\frac{d^{j-1} p_k(x)}{dx^{j -1}} \right|_{x=\sqrt{u}} du}
= 2 \int_0^\infty {\frac{d^{j-1} p_k(t)}{dt^{j -1}}\, t\, dt}.
\end{align*}
Computing for $j=1$, $j=2$, $d \geq j \geq 3$, we have that
\begin{align*}
\int_0^\infty {\MM_1(D_u)du} &= 2\int_0^\infty {p_k(t)t\, dt} = 2 \sqrt{2}\, \frac{\Gamma(\frac{k+1}{2})} {\Gamma(\frac{k}{2})},\\
\int_0^\infty {\MM_2(D_u)du} &= 2\int_0^\infty {p_k'(t)t\, dt} = 2,\\
\mbox{and integration by parts yields}\\
\int_0^\infty {\MM_j(D_u)du} &= 2 \param{ \left. \frac{d^{j-2} p_k(t)}{dt^{j - 2}} t \right|_0^\infty - \left. \frac{d^{j-3} p_k(t)}{dt^{j -3}} \right|_0^\infty} = 0.
\end{align*}
Finally, noting that $\mean{g}=k$ completes the proof.
\end{proof}

\subsubsection{The $F$ case}
Let $M$ be a compact $d$-dimensional manifold and let $f: M \to \R^{n+m}$ be a vector valued Gaussian field with $iid$ components,
\[G(x) = \frac{n}{m} \frac{\sum_{i=1}^m {x_i^2}}{\sum_{i=1}^n{x_{m+i}^2}},\]
and $g = G \circ f$.  In this case, it is proved in \cite[Theorem 15.10.3]{adler2007random} that for $j\ge 1$
\[
\MM_j\param{G^{-1}[u,\infty)} = \param{1+\frac{mu}{n}}^{-\frac{m+n-2}{2}}\sum_{l=0}^{\floor{\frac{j-1}{2}}}\sum_{i=0}^{j-2l-1} C_{m,n,j,l,i} \param{\frac{mu}{n}}^ {\frac{m-j}{2} + i+l}
\]
for a set of constants $C_{m,n,j,l,i}$.

Using basic calculus we can show that for $n> j+2$ and for all $m$, the integral $\int_0^\infty {\MM_j\param{G^{-1}[u,\infty)}du}$ converges. This can be used to compute the expected \textit{lower} Euler integral $\int_M g \ldc$ rather than the expected upper integral that we have computed so far.
Thus, we can conclude that for $n>d+2$ the expected lower Euler integral is finite. For each $n,m$ it is possible to compute the exact value, but no general formula is known. In order to compute the upper Euler integral, we need to compute $\MM_j(G^{-1}(-\infty,u])$. We note that this is feasible, but technically too complicated to be pursued here.

\section{Weighted Sum of Critical Values}\label{s:signed_sum}

Taking $G(x) = H_1(x) = x$ in Theorem~\ref{thm:mean_ei_real}
and using Proposition~\ref{prop:crit_val} yields
the following compact formula.
\begin{theorem}\label{thm:expected_ei_grf}
Let $f: M \to \R$ be a Gaussian random field satisfying the conditions
of the GKF, then
\begin{equation}\label{eq:crazy}
\mean{\int_M{f \udc}} = \mean{\sum_{v\in \CV(f)} \Delta_\chi(f, v) v} = -\frac{\LL_1(M)}{\sqrt{2 \pi}}
\end{equation}
where $\CV(f)$ is set of critical values of $f$ and $\Delta_\chi(f, v)$ is the change
in the {\EC} of $f^{-1}(-\infty, u]$ as $u$ passes through $v$ from below.  In the case that
$M$ is a closed manifold, then
\begin{equation}\label{eq:crazy2}
	\mean{\sum_{p \in \C(f)} (-1)^{\mu(p)} f(p)} = -\frac{\LL_1(M)}{\sqrt{2 \pi}}.
\end{equation}
where $\C(f)$ is the set of critical points of $f$, and $\mu(p) := \mu(p,f)$ is the Morse index of the critical point $p$.
\end{theorem}

In the case that $M$ is a closed even dimensional manifold, $\LL_1(M) = 0$ so \eqref{eq:crazy2} states that
\[\mean{\sum_{p \in \C(f)} (-1)^{\mu(p)} f(p)} = 0.\]

This fact has the following alternative proof, namely:
\begin{align*}
	\mean{\sum_{p \in \C(f)} (-1)^{\mu(p,f)} f(p)} &= \mean{\sum_{p \in \C(-f)} (-1)^{\mu(p,-f)} (-f)(p)}\\
	&= -\mean{\sum_{p \in \C(f)} (-1)^{\mu(p,f)} f(p)}.
\end{align*}
The first equality holds because $f$ and $-f$ have the same law.  The second equality holds deterministically, using that $M$ is even dimensional, because the negative of a Morse function is still a Morse function
and critical points of index $\mu$ become critical points of index $d-\mu$.

The thing to note about Theorem~\ref{thm:expected_ei_grf} is that the expected value of a weighted sum of the critical values
scales like $\LL_1(M)$, a $1$-dimensional measure of $M$ and not the volume $\LL_d(M)$, as one might have expected.
Consider the following example:  Let $f: \R^d \to \R$ be a Gaussian random field with covariance function
$C: \R^d \times \R^d \to \R$ given by $C(x,y) = e^{-\frac{\norm{x-y}^2}{2}}$.  This covariance function induces the Euclidian metric on $\R^d$ and Theorem \ref{thm:expected_ei_grf} implies that
\[
	\mean{\int_{[0,L]^d} f \udc} = -\frac{\LL_1([0,L]^d)}{\sqrt{2 \pi}} = - \frac{d}{\sqrt{2 \pi}}\,L.
\]

In comparison to Theorem~\ref{thm:expected_ei_grf}, letting $G(x)=x^d$ and using Theorem~\ref{thm:mean_ei_real} we get that $\mean{\int_M{f^d \udc}}$  depends on the volume $\LL_d(M)$ (as well the other measures).  So while in general the behavior of the critical points and the critical values depends on the volume, when one takes the weighted sum of the critical values a lot of cancellation occurs and the result only depends on a $1$-dimensional measure.

The result in Theorem~\ref{thm:expected_ei_grf} can be generalized to the case where we consider only critical values below some level $a$.  Observe that taking $a\to\infty$ in the theorem below recovers the result in Theorem~\ref{thm:expected_ei_grf}.
\begin{theorem}\label{thm:mean_signed_sum_level} Let $f: M \to \R$ be a Gaussian random field satisfying the conditions
of the GKF. Then
\begin{equation}\label{eq:crazy3}
\begin{split}
&\mean{\sum_{\stackrel{v\in \CV(f)}{v< a}} \Delta_\chi(f, v) v} =
 -\varphi(a)\LL_0(M) \\
&\mbox{\quad\quad\quad\quad\quad\quad} -\varphi(a)\sum_{j=1}^d{(2\pi)^{-j/2}\LL_j(M)\param{H_{j-2}(-a)- a H_{j-1}(-a)}}.
\end{split}
\end{equation}
In the case that $M$ is a closed manifold, then the left hand side above can be replaced with
$\mean{\sum_{p \in \C(f):\, f(p)<a} (-1)^{\mu(p)} f(p)}$.
\end{theorem}

\begin{proof}
We need to investigate the Euler integral of $G_a(f)$, where recall that
$G_a(x) = x \ind_{(-\infty, a)}(x) + a \ind_{[a, \infty)}(x)$ is the cutoff function from
Corollary~\ref{cor:alter_sum_Ga}.
According to Corollary \ref{cor:alter_sum_Ga},
\[
	\mean{\sum_{\stackrel{v\in \CV(f)}{v< a}} \Delta_\chi(f, v) v}
	= \mean{\int_M G_a(f) \udc} + a\mean{\chi(f \le a)} - a\,\chi(M).
\]
The first term on the right hand side is computed in Lemma \ref{lem:mean_ei_Ga} and the second term is given by the GKF (Theorem \ref{thm:gkf}).
\end{proof}

\begin{lemma}\label{lem:mean_ei_Ga}
Let $f: M \to \R$ be a Gaussian random field, satisfying the GKF conditions. Then,
\[
\begin{split}
\mean{\int_M {G_a(f) \udc}} &= \chi(M)\param{a - a\Phi(a)-\varphi(a)}\\
 &- \varphi(a)\sum_{j=1}^d{(2\pi)^{-j/2}\LL_j(M) H_{j-2}(-a) }
\end{split}
\]
\end{lemma}
\begin{proof}
We will apply the result from Theorem~\ref{thm:mean_ei_general},
so we need to compute $\mean{G_a(f)}$ and
$\int_\R \MM_j(G_a^{-1}(-\infty, u])du$.
Note that
\[
    G_a^{-1}(-\infty,u] = \begin{cases} (-\infty,u] & u < a,\\
     \R & u \ge a.
    \end{cases}
\]
For $j\ge 1$, by \eqref{eq:hermite} we know that $\MM_j((-\infty,u]) = (-1)^{j-1} H_{j-1}(u)\varphi(u)$ and
$\MM_j(\R) = 0$. Therefore if $j\ge 1$, then using that $H_k(-x) = (-1)^k H_k(x)$
\begin{align*}
\int_{\R}{\MM_j(D_u^a) du} &= \int_{-\infty}^a{(-1)^{j-1}H_{j-1}(u)\varphi(u)\, du} \\
&=  \int_{-a}^{\infty}{H_{j-1}(u)\varphi(u)\, du}\\
&= H_{j-2}(-a)\varphi(a),
\end{align*}
where the last transition is due to integration by parts.
We also have that
\[
    \mean{G_a(f)} = \int_{-\infty}^a {x\varphi(x) dx} + \int_a^{\infty} a\varphi(x)dx =  a-a\Phi(a)-\varphi(a).
\]
Thus, by Theorem~\ref{thm:mean_ei_general}, we are done.
\end{proof}

\section{The Persistent Homology of Gaussian Random Fields}

In this section we will give the connection between the Euler integral of a function and its persistent homology.  This will allow us to interpret our computation of the expected Euler integral for Gaussian random fields as a computation on the expected value of a quantitative measure of a Gaussian random field's persistent homology.
We will start off by giving a brief sketch of persistent homology. For more details and further references see
\cite{carlsson2009topology, edelsbrunner2008persistent, ghrist2008barcodes}.

\subsection{Persistent Homology and the Euler Characteristic}

Persistent homology is a way of tracking how the homology of a sequence of spaces changes. For simplicity, in what follows all homology will be with rational coefficients.
Given a filtration of spaces $\mathcal{X} = \{X_u\}_u$ such that $X_s \subset X_t$ if $s < t$, the persistent homology of $\mathcal{X}$, $PH_*(\mathcal{X})$, consists of families of homology classes that `persist' through time.  Explicitly an element of $PH_k(\mathcal{X})$ is a family of homology classes $\alpha = \{\alpha_t\}$ for $t \in [a, b]$, where $\alpha_t \in H_k(X_t)$ (the $k$-th homology group of $X_t$).   These elements are related by the fact that the map
$H_k(X_s) \to H_k(X_t)$, induced by the inclusion $X_s \subset X_t$, maps $\alpha_s$ to $\alpha_t$.  The \emph{birth time}
$a$ of the element $\alpha$ can be thought of as the first time $\alpha$ appears, which is defined by the condition that $\alpha_a$ is not in the image of $H_k(X_s) \to H_k(X_a)$ for all $s < a$.
The \emph{death time} $b$ of the element $\alpha$ is the moment that $\alpha_t$ becomes equivalent to something that existed before $\alpha$.  Formally we require that $\alpha_t$ is not in the image of $H_k(X_s) \to H_k(X_t)$ for all $s < a$ and $t < b$, but $\alpha_b$ is in the image of $H_k(X_s) \to H_k(X_b)$ for all $s< a$.  One must put `tameness' conditions on the filtration so that the birth and death times are defined.

Given a tame function $f: X \to \R$, there is an associated filtration of spaces
$\{f^{-1}(-\infty, u])\}_{u\in\R}$.  One defines the persistent homology of $f$ to be the persistent homology of this filtration, i.e.\! $PH_*(f) = PH_*(\{f^{-1}(-\infty, u])\})$.  The persistent homology of a tame function
$f: X \to \R$ can be seen as a generalization of Morse theory, for if $f$ is a Morse function then the critical values will correspond to birth and death times of elements in the persistent homology $PH_*(f)$.

A graphical way of representing $PH_*(\mathcal{X})$ is via \emph{barcodes}.  Given a persistent homology element $\alpha$, it can be represented by a bar starting at its birth time and ending at its
death time.  To form a barcode, first choose a basis for $PH_*(\mathcal{X})$,
a collection of persistent homology elements such that, for every time $t$, those classes that are alive at time $t$ form a basis for $H_*(X_t)$.  The barcode for $PH_*(\mathcal{X})$ will then be the collection of bars for the chosen basis.

\begin{figure}\label{fig:ph_func}
\centering
\includegraphics[scale=0.32]{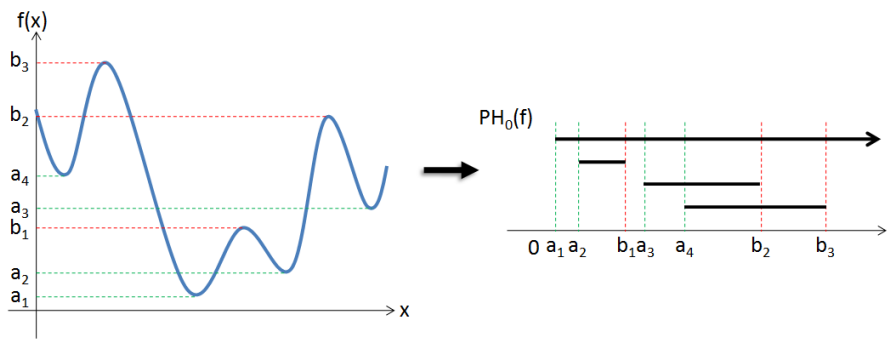}
\caption {The persistent homology of a function $f:\R\to\R$. In this case the sublevel sets are intervals, and therefore the only nonzero homology is $H_0$, representing connected components. The figure to the right presents the barcode of the zeroth persistent homology.  A new $H_0$ element in the persistent homology is born at each local minimum, and a death occurs at each local maximum. When two elements merge, the last one to be born is the first one to die.}
\end{figure}

It turns out that the Euler integral of a tame function $f$ is strongly related to the persistent homology of
$f$.  In light of Proposition~\ref{prop:crit_val}, this is not surprising, since the Euler integral is a measure of how the Euler characteristic of $f^{-1}(-\infty, u]$ changes, while the persistent homology tracks how the homology of $f^{-1}(-\infty, u]$ changes.  To make the relationship precise we need to introduce the following natural extension of the {\EC} to barcodes.

\begin{definition}\label{def:ec_barcode}
	The Euler characteristic of a barcode $B_*$ with a finite number of bars
	and no bars of infinite length is
	\[
		\chi(B_*) = \sum_{b_j \in B_*} (-1)^{\mu(b_j)} \ell(b_j)
	\]
	where $\mu(b_j)$ is the degree of $b_j$ (i.e.\! the homology degree of the class it represents), and $\ell(b_j)$ is the length of $b_j$.
	Equivalently,
	\[
		\chi(B_*) = \int_\R \chi_B(u) du\]
	where $\chi_B(u) = \sum_k (-1)^k \#\set{\textrm{bars of degree $k$ at time $u$}}$ is the signed sum
	of the number of bars at time $u$.
\end{definition}

\begin{proposition} \label{prop:ec_barcode}
	Let $f: X \to \R$ be a tame function and let $PH_*(f,a)$ be the persistent homology of
	$f$ in the range $(-\infty, a]$.  Then
	\[
		\chi(PH_*(f,f_{\max})) = f_{\max}\,\chi(X)- \int_X f \udc,
	\]
	and, in general,
	\[
		\chi(PH_*(f,a)) = a\,\chi(X) - \int_X (G_a \circ f) \udc.
	\]
\end{proposition}
\begin{proof}
Observe that
\begin{align*}
	\int_X f \udc
	&= \int_0^\infty\param{\chi(X)- \chi(f \leq u)} du - \int_{-\infty}^0 \chi(f \leq u)\, du\\
	&= f_{\max}\chi(X) - \int_{-\infty}^{f_{\max}} \chi(f \leq u)\, du.
\end{align*}
Using this equality and $\chi(f \leq u) = \chi_{PH_*(f)}(u)$,
by Definition~\ref{def:ec_barcode} we have
\[
	\chi(PH_*(f,f_{\max})) = \int_{-\infty}^{f_{\max}}{\chi(f \leq u)du}
	= f_{\max}\chi(X) - \int_X f \udc.
\]

As for the second claim, first suppose that $a \leq f_{\max}$.  Then $a = (G_a \circ f)_{\max}$ and
$PH_*(f,a) = PH_*(G_a \circ f, a)$, so by applying the first claim to $(G_a \circ f)$ we get that
\[
	\chi(PH_*(f,a)) = \chi(PH_*(G_a \circ f, a)) = a\,\chi(X) - \int_X (G_a \circ f) \udc.
\]
If $a > f_{\max}$, then
\begin{align*}
	\chi(PH_*(f,a)) = \int_{-\infty}^{a}{\chi(f \leq u)du}
	&= (a-f_{\max})\chi(X) + \int_{-\infty}^{f_{\max}}  \chi(f \leq u)du\\
	&= a\, \chi(X) + \int_X f \udc.
\end{align*}
However, $a > f_{\max}$ implies that $f = G_a \circ f$, so we are done.
\end{proof}

\subsection{The Expected Euler Characteristic of the Persistent Homology of a Gaussian Random Field}

In light of the the connection between the Euler integral of a function and the Euler characteristic of the function's persistent homology in place, we will now reinterpret our computations about the expected Euler integral of a Gaussian random field. This leads to the following result, which as we described in the introduction, seems to be the first result giving a precise form for the expected value of a quantitative property of the persistent homology of random functions.

\begin{theorem}\label{thm:mean_ec_barcode}
Let $f:M\to \R^k$ be a Gaussian random field satisfying the GKF conditions, $G:\R^k \to \R$ continuous and piecewise $C^2$, and $g=G\circ f$. Then
\[
\begin{split}
\mean{\chi(PH_*(g,g_{\max}))} =&  \chi(M)\param{\mean{g_{\max}}-\mean{g}} \\&+ \sum_{j=1}^d{(2\pi)^{-j/2}\LL_j(M)\int_{\R}{\MM_j(D_u) du}}.
\end{split}
\]
If $f:M \to \R$ is a Gaussian random field satisfying the conditions of the GKF, then
\[
\mean{\chi(PH_*(f,f_{\max}))} =  \mean{f_{\max}}\chi(M) + \frac{\LL_1(M)}{\sqrt{2 \pi}}.
\]
\end{theorem}
\begin{proof}
	By Proposition~\ref{prop:ec_barcode},
	\[
		\mean{\chi(PH_*(g,g_{\max}))} = \mean{g_{\max}}\chi(M) - \mean{\int_M g \udc}.
	\]
	Now use the computation of $\mean{\int_M g \udc}$ from Theorem \ref{thm:mean_ei_general}.
\end{proof}

One drawback of this result is that it requires knowledge of $\mean{g_{\max}}$, which is usually unavailable.  However, for real Gaussian random fields, there is a way to circumvent this problem.
\begin{theorem}\label{prop:mean_ec_barcode_bdd}
Let $f: M \to \R$ be a Gaussian random field satisfying the GKF conditions. Then
for any $a \in \R$,
\[
\begin{split}
\mean{\chi(PH_*(f,a))} &= \chi(M)\param{\varphi(a) + a \Phi(a)}\\
&+ \varphi(a) \sum_{j=1}^d{(2\pi)^{-j/2}\LL_j(M) H_{j-2}(-a)}.
\end{split}
\]
\end{theorem}
\begin{proof}
	Again, Proposition~\ref{prop:ec_barcode} gives that
	\[
		\mean{\chi(PH_*(f,a))} = a\chi(M) - \int_M (G_a\circ f) \udc
	\]
	and Lemma \ref{lem:mean_ei_Ga} computes $\int_M (G_a \circ f) \udc$.
\end{proof}

\section{An Application}
An interesting application of the {\EI} is suggested in \cite{baryshnikov2009target}.
Suppose that an unknown number of targets are located in a space $X$, and each target $\alpha$ is represented by its support $U_{\alpha}\subset X$. Suppose also that the space $X$ is covered with sensors, reporting only the number of targets each sensor sees (i.e.\! no identification). Let $h:X \to \Z$ be the \textit{sensor field}, i.e.\! \[h(x) =\#\set{\textrm{targets activating the sensor located at $x$}}.\] The following theorem states how to combine the readings from all the sensors and get the exact number of targets.
\begin{theorem}[\cite{baryshnikov2009target}] \label{thm:target_count}
If all the target supports $U_{\alpha}$ satisfy $\chi(U_{\alpha}) = w$ for some $w\ne 0$, then
\[
\#\set{targets} = \frac{1}{w}\int_X h \udc.
\]
\end{theorem}

Note that we do not need to assume anything about the targets other than they all have the same {\EC}. For example, we need not assume that they are all convex or even have the same number of connected components. On the other hand, the theorem assumes an ideal sensor field, in the sense that the entire (generally continuous) space $X$ is covered with extremely accurate sensors (the range of each sensor is a single point in $X$). In \cite{baryshnikov2009euler} more realizable models using the lower/upper \EI are discussed.

Using the results from Section \ref{sec:ei_grf} we can extend the setup above to the case where the readings from the sensors are contaminated by a Gaussian (or Gaussian related) noise $f(x)$. We will use the following proposition.

\begin{proposition}
Let $h,f: X \to \R$ be tame functions and suppose that $h(X)$ is discrete, then
\[
\int_X (h+f) \udc = \int_X h \udc + \int_X f \udc.
\]
\end{proposition}
\begin{proof}
Let $h(x) = \sum_{i=1}^n{a_i \ind_{A_i}(x)}$, where the $A_i$ are disjoint. Then by the additivity of the {\EC}
we have that
\begin{equation}\label{eq:ei_cover_sum}
\int_X (h+f) \udc = \sum_{i=1}^n  \int_{A_i} (h+f) \udc.
\end{equation}
Next,
\[
	\int_{A_i} (h+f) \udc = \int_{A_i} (a_i+f)\udc = a_i \chi(A_i) + \int_{A_i} f \udc,
\]
where the last equality follows from Proposition~\ref{prop:crit_val}, since every critical value is changed by $a_i$.
Applying this to \eqref{eq:ei_cover_sum} completes the proof.
\end{proof}

Returning to the target enumeration problem, we have a deterministic signal
$x= \int_X h \udc$,
observed via a noisy measurement $Y=\int_X (h+f)\udc$.  By the above proposition we have that
\[Y=\int_X (h+f) \udc = \int_X h \udc + \int_X f \udc=x+N,\] so we have the classical parameter estimation with additive noise model.
If $f(x)$ is a Gaussian or Gaussian related random field satisfying the conditions in Theorem \ref{thm:mean_ei_real}, then we can use the estimator $\hat{x} = Y- \mean{N}$.  Further investigating the properties of the {\EI} might lead to useful estimation techniques for this model.

\section*{Acknowledgments}
The authors would like to thank their supervisors Robert Adler and Shmuel Weinberger, respectively, for their advice and useful comments on this paper.
Thanks are also due to Jonathan Taylor and Yulyi Baryshnikov for helpful discussions during the
Topological Complexity of Random Sets Workshop at AIM, where this work commenced.
\bibliographystyle{plain}
\bibliography{full}
\end{document}